\documentclass{amsart}
\usepackage{appendix} 
\usepackage[colorlinks]{hyperref}        
\usepackage{amsmath, amssymb, mathrsfs, amsfonts,  amsthm} 
\usepackage[all]{xy}

\newcommand{\ham}{\textrm {Ham} (S ^{2}, \omega)}
\newcommand{\diff}{\textrm {Diff} ^{}  (S ^{2})}
\DeclareMathOperator {\LinConfSymp} {LinConfSymp}

\DeclareMathOperator{\Diff}{Diff}

\DeclareMathOperator {\id} {id}
\newtheorem {theorem} {Theorem} [section]

\newtheorem{lemma}[theorem] {Lemma} 
\newtheorem {claim} [theorem] {Claim}
\newtheorem {question}  [theorem] {Question}

\newtheorem {remark} [theorem] {Remark}
\numberwithin {equation} {section}

\begin{document} 
\author {Yasha Savelyev} 
\address{CUICBAS, University of Colima}
\email {yasha.savelyev@gmail.com}
\title {Floer theory and topology of $\Diff (S ^{2} )$}
\begin {abstract} We say that a fixed point of a diffeomorphism is
non-degenerate if 1 is
  not an eigenvalue of the linearization at the fixed point. 
We use pseudo-holomorphic curves techniques to prove the following:
the inclusion map $i: \text{Diff} ^{1} (S ^{2} ) \to \diff
   $ vanishes on all homotopy groups, where $\text{Diff} ^{1} (S
   ^{2} ) \subset  \Diff  (S
 ^{2} )$  denotes the  
space of orientation preserving diffeomorphisms  of $S ^{2} $ with
a prescribed
non-degenerate fixed point. This complements the classical results of
Smale and Eels and Earl.
\end {abstract}
\maketitle
\section {Introduction}
The group of orientation preserving diffeomorphisms of $S
^{2} $ has been shown to deformation retract to $SO (3)$ in 
the classical work of Smale \cite{citeSmaleDiffs}. However it seems
we are still very far from completely understanding some finer aspects
of its topology. Consider for example the following elementary
question. 
\begin{question}
  Let $\Diff^{1}  (S ^{2} )$  denote the space of orientation
  preserving diffeomorphisms of $S ^{2} $, with  one prescribed
  non-degenerate fixed point (with possibly other fixed points). Here non-degenerate means that 1 is
  not an eigenvalue of the linearization at the fixed point. Is this
  space contractable?
\end{question}
We expect that the answer is yes, since $\Diff^{1}  (S ^{2} ) \cap SO
(3)$ is just the space $SO ^{1}  (3)$ of rotations about a fixed axis, not containing the identity
map, and so forms a contractible space. Moreover Earle and Eells
\cite{citeEarleEellsAfiberBundleDescription} prove a somewhat related fact that the space of orientation preserving
diffeomorphisms of $S ^{2} $ fixing three points is contractible.

The intuition
behind the argument of Smale or Earle and Eells is indicative
that there should be a deformation retraction of $\Diff^{1}  (S ^{2} ) $ onto $SO ^{1}  (3)$.   One tricky point in directly applying Smale's
argument (for instance) is that in constructing the deformation
retraction  one must be careful to keep our fixed point fixed and
non-degenerate.  If we allow the deformation retraction to take
place in the ambient space $\Diff  (S ^{2} ) $, of orientation
preserving diffeomorphisms,  some version of Smale's argument
may (very likely) give  a null-homotopy of the inclusion $\Diff ^{1} (S
^{2} ) \to \Diff (S ^{2} ) $, although this is not known at the
moment. 
We show using very different techniques the following:
\begin{theorem} \label{thm:main2} The inclusion map $i: \text{Diff}
   ^{1} (S ^{2} ) \to \diff
   $ vanishes on all homotopy groups.
\end{theorem}
 Our proof involves elliptic, or more precisely pseudo-holomorphic curve techniques in
symplectic geometry, and  is morally closely connected to Floer theory.    We
certainly hope that our techniques  can be sharpened to give that
$\text{Diff} ^{1} (S
^{2}, \omega ) $ is itself contractible. 

Part of the interest in
developing the main technique of this paper is that some version of it may be
applicable in some other contexts, particularly the context of
contactomorphism groups of contact 3-folds.


\begin{remark}
A technically identical argument but dealing with Hamiltonian fibrations over
discs, rather than $\mathbb{CP} ^{1} $, gives an analogue of
the above theorem  for
the space $\text{Lag}(S ^{2} ) $ of
Lagrangians submanifolds in $S ^{2} $ Hamiltonian isotopic to equator,
i.e.
simple closed smooth unparametrized
loops  partitioning $S ^{2} $ into a pair of regions with equal
area. The basic ingredient for this is that a Hamiltonian loop of
Lagrangian submanifolds, that is a loop obtained by a Hamiltonian flow,
gives rise to a sub-fibration of $M \times D ^{2} $ over the boundary
$\partial D ^{2} $  with fiber $L _{\theta}  \subset M _{\theta} $ over
$\theta$ the element of the loop at $\theta$. We may then analogously study moduli
spaces of holomorphic disks with boundary on this Lagrangian
sub-fibration, for suitable almost complex structures, see for instance
\cite{citeAkveldSalamonLoopsofLagrangiansubmanifoldsandpseudoholomorphicdiscs.}
for a similar story, or
\cite{citeHuLalondeArelativeSeidelmorphismandtheAlbersmap.} for the
monotone case.
\begin{claim}
The inclusion map $i: \text{Lag} ^{1} (S ^{2} ) \to \text{Lag} (S
^{2} )$ vanishes on all homotopy groups.  Where the superscript
$1$ means that we take the subspace of
  those Lagrangians transversally intersecting the standard equator at
  a prescribed point, (there may be other intersections).
  \end{claim}
\end{remark}
\section {Proof of Theorem \ref{thm:main2}}
Let $\Omega  \text{Diff} ^{1}  (S ^{2})$
denote the space of based smooth loops in $\text{Diff} ^{1}  (S ^{2}, \omega )$,   which are
constant near end points.
We show that $i$ is null-homotopic as follows. Associated to each
$p \in \Omega  \text{Diff} ^{1}  (S ^{2}) $
there is a structure group $\diff$ $S ^{2}  $ fibration $X _{p} $ over $\mathbb{CP}
^{1}  $, with $\diff$ the group of orientation preserving
diffeomorphisms,  which is 
formed via the clutching construction. It will be more illuminating to
think of it as a conformal symplectic fibration, that is a fibration
whose transition maps are conformal symplectic.

 The fact that we are pushing $g$ from $\Omega \text{Diff} ^{1}  (S ^{2})$, specifically the
non degeneracy property of the associated diffeomorphisms is used
to construct a very special almost complex structure on $X _{s}= X
_{i_*g (s)}$ parametrically with $s$ and using that a natural
foliation of $X _{s} $ by holomorphic
sections.
This foliation in turn
determines a natural smooth trivialization of $X _{s} $. The construction of this
foliation uses in particular classical
positivity of intersections ideas of Gromov-McDuff, particularly following
\cite{citeGromovPseudoholomorphiccurvesinsymplecticmanifolds.}.
Although our techniques are based on the theory of closed
pseudo-holomorphic curves, the
underlying idea is closely related to Floer theory. In particular we have in mind some potential
generalizations which use Floer theory more explicitly.
Given all this 
it is a matter of topology to deduce the main result. 

We now give the detailed argument. For each $s$ we
get a conformal symplectic $S ^{2} $
fibration $$\pi _{s}: X _{s} \to \mathbb{CP} ^{1},$$ as
follows 
\begin{equation} \label{eq:clutch}
X _{s} : = S ^{2} \times D ^{2} _{-}   \sqcup _{g (s)}  S ^{2} \times D
^{2} _{+} ,
\end{equation}
$D ^{2} _{\pm}  $ identified with $D ^{2} $, and $X _{s} $
 is the quotient of $S ^{2} \times D ^{2} _{-} \sqcup S ^{2} \times D
 ^{2} _{+}  $ by the equivalence relation:
$$(x, (1, \theta)) \in S ^{2}
 \times \partial D ^{2} _{-}    $$ is equivalent
to $$(g (s)   (\theta) ^{-1}
(x), (1, \theta)) \in S ^{2} \times \partial D ^{2} _{+}   $$ 
using the polar coordinates $(r, \theta)$, $r \in [0,1]$ on $D ^{2} $.
The total space of this $s$-family can be understood as a smooth
conformal symplectic $S ^{2} $ fibration $$\widetilde{\pi}:
\widetilde{X} \to S ^{m} \times S ^{2}.  $$
For each $s$ we have a 
smooth  section $\sigma _{s}  $ of
$X _{s}  $  corresponding to our distinguished non-degenerate fixed point $x _{fix}
$, in our coordinates over $D ^{2} _{\pm}  $ this section is just $z \mapsto x
_{fix} $.
We shall call these  
 \emph{spectral sections}. 

\begin{lemma} \label{lemma:verticalChernnumber0} 
   The spectral sections $\sigma _{s} $ have vertical Chern number 0.
\end{lemma}
\begin{proof}  
 By construction the  restriction $N _{s} $ of the vertical bundle of $X _{s} $ to $\sigma
_{s} $ is obtained by the clutching construction $$N _{s} \simeq \mathbb{R} ^{2}
 \times D ^{2} _{-} \sqcup _{g _{s,*}}  \mathbb{R} ^{2} \times D ^{2}
 _{+},    
 $$ where $g _{s,*}$  denotes the linearization of the loop $f _{s} $
 at $x _{fix} $. The Conley-Zehnder index obviously extends to the
 linear conformal symplectic group $\LinConfSymp (\mathbb{R} ^{2n} )$, 
 and  $ CZ (g _{s, *} )=0$ where $CZ$ is
 the Conley-Zehnder index,   and $CZ (g _{s,*} )=0$
 directly follows by
 the assumption that $x _{fix} $ is a non-degenerate fixed point for
 $g (s) (\theta)$ for each $\theta$, and by the construction of the
 Conley-Zehnder index in
 \cite{citeRobbinSalamonTheMaslovindexforpaths.}, (we have a loop  with no crossing points
 with the Maslov cycle).  On the other
 hand $CZ (g _{s,*} )$ is exactly the Chern number of $N _{s} $, (for
 suitable normalizations.) See
 for example \cite{citeMcDuffSalamonIntroductiontosymplectictopology},
 actually it is obvious that this is true up to a non-zero multiple,
 since both the Chern number and the Conley-Zehnder index determine
 injections $\pi _{1} (\LinConfSymp (\mathbb{R} ^{2n})) \to \mathbb{Z} $,
 which is all that is necessary here. 
 \end{proof}

Next we construct a smooth family  $\{\mathcal{A}_{s} \}$  of
$\diff$ connections on
$\{X _{s} \}$, as follows. Let $\mathcal{F} \to S ^{k} $ be a
fibration whose fiber $\mathcal{F} _{s} $ over $s$ is the space of
$\diff$ connections on
$X _{s} $, for which the spectral section $\sigma _{s} $ is flat. Each
$\mathcal{F}_{s} $ is a non-empty, affine space and hence contractible. In
particular we may choose a smooth section $\{\mathcal{A}_{s} \}$ and this is our
family. 

Fix a smooth family 
$\{j _{z,s} \}$, $z \in \mathbb{CP} ^{1}, s \in S ^{m}  $ of
$\{\omega _{z,s} \} $-compatible almost
complex structures on the vertical tangent bundles $\{T ^{vert}|
_{z,s} \widetilde{X}   \}$. 
A smooth $\diff$ connection $\mathcal{A}$ on $X _{s} $,  gives rise to an
almost complex structure $J _{\mathcal{A}} $, by asking that
horizontal spaces be $J _{\mathcal{A}} $-invariant, that on $T
^{vert}| _{z} X _{s}   $ $J _{\mathcal{A}} $ coincides with $j _{z,s}
$, and so that the projection map $\pi _{s} $ is $J _{\mathcal{A}}
$-holomorphic.
We let $\{J
_{s} \} $ denote the family of  almost complex
structures on $ \{X _{s} \}$ corresponding to $\{\mathcal{A} _{s} \}$.
\begin{lemma} \label{lemma:compact}
   The moduli space
   $\overline{\mathcal{M}} ( J _{s}, A:=[\sigma _{s} ]  )$
of  
stable $J _{s}   $-holomorphic sections of $X _{s} $, in the class   $A:=[\sigma
_{s} ]$,
contains no nodal curves. 
\end{lemma}
By a \emph{stable holomorphic section} we mean a stable $J _{s} $-holomorphic
map $\sigma$ into $ X _{s} $, in the classical sense,
\cite{citeMcDuffSalamon$J$--holomorphiccurvesandsymplectictopology}
with domain an unmarked nodal Riemann sphere, one
of those components is called \emph{principal}. The restriction $\sigma
_{princ} $ of $\sigma$, to the principal component is a $J _{s} $-holomorphic section, i.e. we have a commutative diagram:
\begin{equation*}
\xymatrix {& X _{s} \ar [d] ^{\pi _{s} }  \\ 
\mathbb{CP}^{1} \ar [ur] ^{\sigma _{princ} } \ar [r] ^{\id}  &   \mathbb{CP}.}
\end{equation*}
All the other components of $\sigma$ are vertical, that is they are
$J_s$-holomorphic maps into the fibers of $X _{s} $. 
\begin{proof}
By construction we have a smooth, embedded $\mathcal{A} _{s} $-flat and hence
$J _{s} $-holomorphic section
$\sigma _{s}  $ of $X _{s} $.
 If there are stable sections with more than one component in our moduli space
   $\overline{\mathcal{M}} ( J _{s}, A   )$ 
 then taking the
 principal component of the corresponding nodal section, we would have a smooth
 holomorphic section $\sigma'$ in class $[\sigma _{s} ] + A$ for some spherical fiber class $A$ with $c _{1} (A)<0$. 

 Consequently since $[\sigma _{s} ] \cap [\sigma _{s} ] =0 $ by Lemma
 \ref{lemma:verticalChernnumber0}, $\sigma'$ would have negative intersection number with
 $\sigma _{s} 
  $,
  which is impossible by positivity of intersections as $\sigma _{s}   $ is
  embedded, see \cite [Section
  2.6]{citeMcDuffSalamon$J$--holomorphiccurvesandsymplectictopology}.
  \end{proof}
\begin{lemma}
   $\overline{\mathcal{M}} ( J _{s}, A   )$  is also regular, (in
  particular has the expected dimension). 
\end{lemma}
\begin{proof}
   By automatic transversality the top stratum is regular, as since
   all elements (in the top stratum) have normal bundle with vanishing
   Chern number, the associated real linear CR
  operator is automatically surjective, see for example
  \cite[Appendix
  C]{citeMcDuffSalamon$J$--holomorphiccurvesandsymplectictopology}. 
  But we just showed that there are no nodal curves and hence no other
  strata. 
  %
  \end{proof}

Since we have no nodal curves, all elements of
$\overline{\mathcal{M}} ( J _{s} , A)$ are in particular embedded, and
this moduli space is non-empty as it has the
element $\sigma _{s} $, then by positivity of
intersections, and by regularity 
\begin{equation*}
   ev _{s}: \overline{\mathcal{M}} ( J _{s} , A) \to \pi _{s} ^{-1} (z _{0} ) \simeq S ^{2},
\end{equation*} 
taking a section to its value over $0 \in \mathbb{C} $ 
must be a smooth degree 1 map. But then again by positivity of
intersections, 
we get that the sections of $\overline{\mathcal{M}} ( J _{s} , A) $ induce a
canonical smooth foliation of $X _{s} $,
with all leaves diffeomorphic to $\mathbb{CP}^{1}$, and as they are sections we
get a canonical smooth trivialization $$X _{s} \xrightarrow{tr _{s} } S ^{2} \times
\mathbb{CP}^{1},$$  by identifying all fibers with the fiber over $0 \in
\mathbb{CP}^{1}  $ using the leaves of the foliation. In particular
each $X _{s} $ is trivial as a smooth fibration.
 So for the rest of the argument we
assume that $m>0$.

The total space $\widetilde{X} $ of the family $\{X _{s} \}$ can be understood as a bundle $P _{g}  \to S ^{m}    $  with fiber over $s$: $X
_{s} \simeq S ^{2} \times \mathbb{CP} ^{1}   $ and
with structure group by the discussion above reducible to the group of
$\diff$ bundle automorphisms
of $S ^{2} \times \mathbb{CP}^{1}$ fixing the fiber over $0$, with $0$
identified with the
origin of $D ^{2} _{-}  $.  In
other words this structure group is: 
\begin{equation*} 
\Omega ^{2} _{\id} \diff,
\end{equation*}
with $\Omega ^{2} $
denoting the based spherical mapping space in the class of the
constant map.
We drop the subscript $\id$ from now on.

By construction $P _{g} $ is equivalent to the pullback by the sequence
\begin{align*}
   S ^{m}   \xrightarrow{g}
\Omega  \text{Diff} ^{1}  (S ^{2}, \omega ) \xrightarrow{i_*} 
   \Omega  \diff,
\end{align*}
of the structure group
$\Omega ^{2} \diff$ bundle $U$ over $\Omega \diff $, whose fiber 
over a loop $\gamma$
is the fibration $X _{\gamma} $ constructed as in
\eqref{eq:clutch}. 
 It follows by the author's
\cite [Section 7]{citeSavelyevQuantumcharacteristicclassesandtheHofermetric.} that
$U$ is the universal structure group $\Omega ^{2} \diff $
bundle. To note, we state results there for (as a particular case) $\ham$ rather than $\diff$ but the argument 
is the same for $\diff$. Strictly speaking what follows by
\cite{citeSavelyevQuantumcharacteristicclassesandtheHofermetric.} is
that $U$ is universal for the subgroup of $\Omega ^{2} \diff $
corresponding to bundle automorphisms (over $\id$) of $S ^{2} \times \mathbb{CP}
^{1}  $, which are identity bundle maps over contractible
neighborhoods of 
$0, \infty \in \mathbb{CP} ^{1} $, however this subgroup can be easily shown to be homotopy equivalent to $\Omega
 ^{2} \diff$.  

So if $P$ is trivial as   a structure group $\Omega ^{2}
\diff$ bundle the map $i_* \circ g$ is null-homotopic. 
\qed

\section {Acknowledgements}
The paper was completed at ICMAT Madrid,
where the author was a research fellow, with funding from Severo Ochoa grant and
Viktor Ginzburg laboratory. I would like to thank  
Michael Usher for reading
of an early draft and pointing out a gap,
Dusa McDuff for interest, and much help in
clarifying the paper and Viktor Ginzburg  for
interest and
comments. As well as the anonymous referee for very careful reading
and finding some serious issues.
 \bibliographystyle{siam}  
 \bibliography{/root/texmf/bibtex/bib/link} 
\end{document}